\documentclass{lmcsArX} 
\keywords{stochastic SIS epidemic model, explicit numerical scheme, semi-discrete method, strong first order convergence, exponential stability}

\usepackage{hyperref}
\theoremstyle{plain} 

\def\cf{{\em cf.}}

\usepackage{amsmath,amsfonts,amsbsy,latexsym,multicol,fancybox}
\usepackage{mathtools}
\DeclarePairedDelimiter\ceil{\lceil}{\rceil}

\usepackage{graphicx}
\usepackage{epstopdf}
\usepackage{bm}

\newtheorem{theorem}{Theorem}
\newtheorem{assumption}{Assumption}
\newtheorem{definition}{Definition}
\newtheorem{lemma}{Lemma}
\newtheorem{remark}{Remark}

\newcommand{\bass}{\begin{assumption}}\newcommand{\eass}{\end{assumption}}
\newcommand{\bde}{\begin{definition}} \newcommand{\ede}{\end{definition}}
\newcommand{\ble}{\begin{lemma}} \newcommand{\ele}{\end{lemma}}
\newcommand{\bth}{\begin{theorem}} \newcommand{\ethe}{\end{theorem}}
\newcommand{\bre}{\begin{remark}} \newcommand{\ere}{\end{remark}}

\newcommand{\bpf}{\begin{proof}}\newcommand{\epf}{\end{proof}}
\newcommand{\barr}{\begin{array}}\newcommand{\earr}{\end{array}}
\newcommand{\beao}{\begin{eqnarray*}}\newcommand{\eeao}{\end{eqnarray*}\noindent}
\newcommand{\beam}{\begin{eqnarray}}\newcommand{\eeam}{\end{eqnarray}\noindent}
\newcommand{\beqq}{\begin{equation}}\newcommand{\eeqq}{\end{equation}\noindent}

\newcommand{\wt}{\widetilde}

\newcommand{\tto}{t\to\infty}

\newcommand{\be}{\beta}
\newcommand{\ga}{\gamma} 

\newcommand{\Del}{\Delta}

\newcommand{\w}{\omega} \newcommand{\W}{\Omega}

\newcommand{\bfE}{{\mathbb E}} 
\newcommand{\bbf}{{\mathcal F}}

\newcommand{\bbl}{{\mathcal L}}

\newcommand{\bfP}{{\mathbb P}}

\newcommand{\bbr}{{\mathcal R}} \newcommand{\bbR}{{\mathbb R}}

\usepackage{lineno}
\usepackage{doi}

\usepackage{tikz}
\usetikzlibrary{cd}
\usetikzlibrary{matrix}
\usetikzlibrary{automata}
\usetikzlibrary{arrows}
\usepackage{rotating}
\usepackage{subcaption}

\usepackage{color}

\definecolor{pblue}{rgb}{0.13,0.13,1}
\definecolor{pgreen}{rgb}{0,0.5,0}
\definecolor{pred}{rgb}{0.9,0,0}
\definecolor{pgrey}{rgb}{0.46,0.45,0.48}

\begin{document}

\title[Explicit scheme for the stochastic SIS epidemic model]{Domain preserving and strongly converging explicit scheme for the stochastic SIS epidemic model}

\author[Y.~Kiouvrekis]{Yiannis Kiouvrekis\lmcsorcid{0000-0001-6805-3203}}[a,c]

\address{Mathematics, Computer Science, Artificial Intelligence Laboratory, Department of Public and One Health, 
University of Thessaly}	
\email{kiouvrekis.y@uth.gr}  

\author[I.S.~Stamatiou]{Ioannis S. Stamatiou\lmcsorcid{0000-0002-8215-9634}}[a,b]
\address{University of West Attica, Athens, Greece}
\email{istamatiou@uniwa.gr}

\address{University of Nicosia, Business School, Cyprus}




\begin{abstract}
  \noindent In this article, we construct a numerical method for a stochastic version of the Susceptible–Infected–Susceptible (SIS) epidemic model, expressed by a suitable stochastic differential equation (SDE), by using the semi-discrete method to a suitable transformed process. We prove the strong convergence of the proposed method, with order $1,$ and examine its stability properties.  Since SDEs generally lack analytical solutions, numerical techniques are commonly employed. Hence, the research will seek numerical solutions for existing stochastic models by constructing suitable numerical schemes and comparing them with other schemes. The objective is to achieve a qualitative and efficient approach to solving the equations. Additionally, for models that have not yet been proposed for stochastic modeling using SDEs, the research will formulate them appropriately, conduct theoretical analysis of the model properties, and subsequently solve the corresponding SDEs.
\end{abstract}

\maketitle


\section*{Introduction}\label{S:one}

The research on mathematical modeling of epidemics initially focused on deterministic models, which significantly contributed to understanding epidemic behavior \cite{BRAUER2017113}. Deterministic models were able to capture various disease characteristics, such as permanent immunity and sexually transmitted or bacterial diseases without permanent immunity \cite{c18ef097-75bd-3f59-a5d7-256a8f9b744e}. However, deterministic models have limitations in adequately describing real-world scenarios due to the influence of uncertain circumstances on model parameters. To address this limitation, stochastic models were introduced, incorporating noise factors to better suit epidemiological problems \cite{Greenwood2009}. This paper aims to explore the modeling of population dynamics using stochastic differential equations (SDEs). The focus is not only on formulating the initial model as an appropriate SDE using independent Brownian motions but also on developing numerical approaches to solve the equations. \\
SDEs, as we mentioned, play a prominent role in several areas as finance and epidemiology but in comparison with ODE the general solution theory is more mathematical complicated \cite{sauer}. The same situation exists in the case of numerical solutions of SDEs. Methods for the numerical solution of SDEs are based on similar techniques like the Euler Maruyama method (EM)\cite{maruyama}, which is the analogue of the Euler method for ordinary differential equations.
It is well known that if the SDE’s coefficients are globally Lipschitz continuous, then the Euler approximation process will convergence in the strong and numerically weak sense to the exact solution of the SDE, but if the coefficients of the SDE are not globally Lipschitz continuous the Euler approximation does not converge \cite{hutz}.\\
The initial reference is the research paper \cite{gray_elal:2011} where the authors have expanded upon the classical susceptible-infected-susceptible (SIS) epidemic model. The deterministic model which also describes the vital dynamics of the population \cite{c18ef097-75bd-3f59-a5d7-256a8f9b744e} is the following: 
\begin{equation} \label{eq1}
\begin{split}
\dfrac{dS(t)}{dt} & = -\frac{\beta}{K} S(t)I(t) + (b + \gamma) I(t) \\
\dfrac{d I(t)}{dt} & = \frac{\beta}{K} S(t)I(t) - (b+\gamma)I(t)
\end{split}
\end{equation}
where $\beta>0$ is the contact rate, $\gamma>0$ the recovery rate, $b\geq0$ the birth rate for a population of size $K$ with initial conditions $S(0) + I(0) = K$ with $S(0)>0$ and $I(0)>0.$ 

The authors in \cite{gray_elal:2011} have taken this model from its original deterministic framework and transformed it into a stochastic framework (\ref{eq2}). To accomplish this, they have employed stochastic differential equations to describe the dynamics of the number of infectious individuals, represented by $I(t)$.
\begin{equation} \label{eq2}
\begin{split}
dS(t) & =\left[-\frac{\beta}{K} S(t)I(t) + (b + \gamma) I(t)\right]dt-\sigma S(t)I(t)dW(t)\\
dI(t) & = \left[\frac{\beta}{K} S(t)I(t) - (b+\gamma)I(t)\right]dt + \sigma S(t)I(t)dW(t)
\end{split}
\end{equation}
where $W(t,\w\footnote{we usually omit the dependence on $\w$}):[0,\infty)\times\W\rightarrow\bbR$ is an $1$-dimensional Wiener process adapted to the filtration $\{\bbf_t\}_{t\geq 0},$ see \cite{karatzas_shreve:1988}, \cite{mao:2007} and $\sigma$ belongs to $\mathbb{R}^{+}$. Expressing the SDE in terms of the process $I$ in integral form we get that 
\beqq\label{eq3}
I_t = I_0 + \int_0^t\left(  \eta I_s - \frac{\beta}{K} I^2_s \right) ds + \sigma\int_0^t (K - I_s)I_sdW_s,
\eeqq
where $\eta = \beta - b - \gamma.$ It has been shown in \cite{gray_elal:2011} that for any given initial value $I_0 = I(0)\in(0,K)$ there exists a unique global solution with values in $(0,K)$ in the sense that $\bfP(I_t\in(0,K) \hbox{ for all } t\geq0) = 1.$ Unfortunately, the classical EM method does not preserve the domain $(0,K).$ 

The epidemic models refer to quantities which take values in a certain domain; for the SIS model $I(t)\in[0,K].$ Therefore, the numerical scheme should preserve this domain. The commonly used explicit Euler scheme does not have this property since its increments are conditional Gaussian and thus there is an event of negative values with positive probability. We are interested in the construction of a numerical scheme that preserves the domain of the solution process. Of course other features of the numerical method (apart from the domain preservation) are desirable, such as the strong convergence in the mean square sense to the exact solution of the original SDE (for visualizing stochastic dynamics, simulating scenarios theoretical interest, see \cite{hutzenthaler_jentzen:2015}) and the explicitness of the method (for computational reasons). 

From all of the above, the need to approximate in a qualitative correct way the solution processes in the nonlinear models as in  the SIS model emerges (\ref{eq3}). The idea of domain-preserving numerical methods is a direction where researchers have paid attention, more obvious in the last fifteen years, c.f.  
\cite{dangerfield:2012}, \cite{kahl_et_al:2008}, \cite{kelly_et_al:2018}, \cite{halidias:2015}, \cite{halidias_stamatiou:2023} and \cite{STAMATIOU:2021} and references therein. If in addition we require explicit numerical schemes we could use the semi-discrete method,  see \cite{STAMATIOU:2021} for a review of the method.

The semi-discrete method, originally proposed in \cite{halidias:2012}, has the following properties: 

\begin{itemize}
	\item it is in general explicit in general and consequently does not require a lot of computational time,
	\item it does not explode in non-linear problems,  
	\item it strongly converges to the exact solution of the original SDE, 
	\item domain preservation,
 \item reproduces the stability behavior of the solution process (\cite{halidias_stamatiou:2022}, \cite{STAMATIOU:2021b})
\end{itemize}

The key idea behind the semi-discrete method is freezing on each integration interval of size $\Del,$ parts of the drift and diffusion coefficients of the solution at the endpoints of the subinterval, obtaining explicitly solved SDEs which by construction preserve the domain of the solution process.

For the SIS model (\ref{eq2}) the following domain preserving method has been proposed, \cite{YANG:2021},

\beqq\label{eq:SISGrayYang}
Y_{n} = \frac{Ke^{X_n}}{1 + e^{X_n}}, 
\eeqq
where $X_n$ is produced by the application of the EM scheme to the transformed process, through the Lamperti-type transformation $z_t = \ln \dfrac{I_t}{K - I_t},$ with  

\beqq\label{eq:SIStransformed}
z_{t} = z_0 + \int_0^t F(z_s) ds + \int_0^t \sigma K d W_s, 
\eeqq
where 
$$
F(x) = \eta - (b + \gamma)e^x + \frac{\sigma^2 K^2}{2} - \frac{\sigma^2 K^2}{1 + e^x}
$$
that is, 
\beqq\label{eq:SISGrayYangTransf}
X_{n+1} = X_n + F(X_n) \Del + \sigma K\Del W_n, 
\eeqq
with $X_0 = z(0) = \ln \dfrac{I_0}{K - I_0}.$ 

We propose the following domain preserving scheme applying the semi-discrete method to a different transformation of the original process and then transforming back. In particular, the numerical scheme we propose is produced by 

\beqq\label{eq:SIS_SD}
\hat{Y}_{n} = \frac{K\hat{X}_n}{1 + \hat{X}_n}, 
\eeqq
where $\hat{X}_n$ is produced by the application of the semi-discrete method to the transformed process $\hat{z}_t = \dfrac{I_t}{K - I_t},$ with  

\beqq\label{eq:SIS_SD_transformed}
\hat{z}_{t} = \hat{z}_0 + \int_0^t \hat{F}(\hat{z}_s) ds + \int_0^t \sigma K \hat{z}_s d W_s, 
\eeqq
where 
$$
\hat{F}(x) = \eta x - (b + \gamma)x^2 + \frac{\sigma^2 K^2 x^2}{1 + x} =: x\phi(x)
$$
that is, 
\beqq\label{eq:SIS_SD_Transf}
\hat{X}_{n+1} = \hat{X}_n\exp\{ (\phi(\hat{X}_n) - \frac{\sigma^2 K^2}{2}) \Del + \sigma K\Del W_n\}, 
\eeqq
with $\hat{X}_0 = \hat{z}(0) = \dfrac{I_0}{K - I_0}.$ Exponential strongly converging schemes like (\ref{eq:SIS_SD_Transf}) had been proposed in \cite{halidias_stamatiou:2016} but with no rate of convergence. The proposed numerical scheme (\ref{eq:SIS_SD}) is proven to strongly converge to the solution process with order $1,$ see Theorem~\ref{th1}. Moreover, the numerical scheme possesses another property, reproduces the stability behavior of the  solution process with the cost of an extra mild assumption on the parameters, see Theorem~\ref{th2}. 

In Section~\ref{sec:Main} the setting and the two main results are presented. The proofs are found in Sections~\ref{sec:Con} and \ref{sec:Stab}. Finally, Section~\ref{sec:Num} provides numerical evidence of the validity of the theoretical results. 

\section{Setting, Preliminary and Main results}\label{sec:Main}

Recall the SDE (\ref{eq3}) for the process $I$  which we rewrite as

\beqq\label{eq:SIS_forI}
I_t = I_0 + \int_0^t A(I_s) ds + \int_0^t B(I_s) dW_s,
\eeqq
with 
\beqq\label{eq:SIS_Icoefficients}
A(x) = \eta x - \frac{\beta}{K} x^2 \quad\hbox{ and }\quad  B(x) = \sigma(K - x)x,
\eeqq
with $\eta = \beta - b - \gamma.$ 
The diffusion operator associated with (\ref{eq:SIS_forI}) related to the transformation $V(x) = \dfrac{x}{K - x}$ reads
\beao
LV(x) &=& A(x)V^\prime(x) + \frac{1}{2}B^2(x)V^{\prime\prime}(x)\\
& = & (\eta x - \frac{\beta}{K} x^2)\dfrac{K}{(K - x)^2} + \frac{1}{2}\sigma^2(K - x)^2x^2\dfrac{2K}{(K - x)^3}\\
& = & \eta \dfrac{x}{K - x}\dfrac{K}{K - x}   - \beta\dfrac{x^2}{(K - x)^2} + \sigma^2 K\dfrac{x}{K - x}x\\
& = & \eta V(x)(1+V(x)) - \beta)V^2(x) + \sigma^2 KV(x) \dfrac{KV(x)}{1+V(x)}\\
& = & (\eta - \beta)V^2(x)+ \frac{\sigma^2 K^2 V^2(x)}{1 + V(x)}
\eeao
therefore by application of the It\^{o} formula we get (\ref{eq:SIS_SD_transformed}), rewritten as

\beam
\nonumber\hat{z}_{t} &=& \hat{z}_0 + \int_0^t LV(I_s) ds + \int_0^t B(I_s)V^\prime(I_s) d W_s\\
\label{eq:SIS_SD_transformedAgain}& = & \hat{z}_0 + \int_0^t \hat{F}(\hat{z}_s) ds + \int_0^t G(\hat{z}_s) d W_s 
\eeam
with  
\beqq\label{eq:SIS_Zcoefficients}
\hat{F}(x) = \eta x - (b + \gamma)x^2 + \frac{\sigma^2 K^2 x^2}{1 + x} \quad\hbox{ and } \quad G(x) = \sigma K x
\eeqq

\ble
Process $\hat{z}_{t}$ has finite moment bounds of any order, that is for any $p\in\bbR$ there is a constant $C_p$ such that $$\sup_{t\in[0,T]}\bfE(\hat{z}_{t}^p)\leq C_p$$
\ele
\bpf
Take a $p>0$ and write $\hat{z}_{t} = I_t(K-I_t)^{-1}.$ Raising to the power of $p$ and taking expectations 
$$
\sup_{t\in[0,T]}\bfE \left(\hat{z}_{t}^p\right) \leq \sqrt{\sup_{t\in[0,T]}\bfE (I_t)^{2p}}\sqrt{\sup_{t\in[0,T]}\bfE (K - I_t)^{-2p}} \leq K^p \sqrt{\hat{C}_{2p}},
$$
whereas
$$
\sup_{t\in[0,T]}\bfE\left( \hat{z}_{t}^{-p}\right) \leq \sqrt{\sup_{t\in[0,T]}\bfE (I_t)^{-2p}}\sqrt{\sup_{t\in[0,T]}\bfE (K - I_t)^{2p}} \leq K^p \sqrt{\hat{C}_{2p}},
$$
with 
$$\hat{C}_{p} = \left((I_0)^{-p}\vee(K- I_0)^{-p}\right)\exp\left\{p\left(\eta\vee(2\beta-\eta)\vee \frac{2\beta}{K}\right)T + \frac{p(p+1)}{2}\sigma^2K^2T\right\}$$ 
being the constant as in \cite[Theorem 3.2]{chen_etal:2021}. 

Alternatively, for all positive $x$ we find that
\beao
J(x)&:=&\frac{x\hat{F}(x) + (p-1)G^2(x)/2}{1 + x^2}\\
& = &\frac{x^2}{1 + x^2}\left(\eta - (b+\ga)x + \frac{\sigma^2 K^2 x}{1 + x} + \frac{(p-1)\sigma^2K^2}{2}\right)\\
& \leq&C,
\eeao
for any $p>2,$ where from now and on $C$ refers to a constant varying form line to line.  Therefore, the $p$-th moment of $(\hat{z}_{t})$ is finite, $\bfE (\hat{z}_{t}^{p})< C$ \cf \cite[Theorem 2.4.1]{mao:2007}.  Following a standard procedure, we may show that $\bfE\sup_{t\in[0,T]} (\hat{z}_{t}^{p})<C,$ c.f.\cite[Lemma 4.3]{halidias_stamatiou:2016} for $p>2$ and work in a similar way for the process $(\hat{z}_{t})^{-p}.$ The case where $0<p<2$ is covered by Jensen's inequality for the concave function $x^{p/2}.$
\epf

Now let us discuss about approximation schemes for the solution of (\ref{eq:SIS_SD_transformedAgain}). Using the auxiliary function $f(x,y) = \left(\eta  - (b + \gamma)x + \frac{\sigma^2 K^2 x}{1 + x}\right)y:= \phi(x)y,$ with the property $f(x,x) = \hat{F}(x),$ we write the evolution of a process $(\hat{X}_{t})$ in a subinterval $[t_n,t_{n+1}]$ of length $\Del<1$ as 
\beao
\hat{X}_s &=& \hat{X}_{n} + \int_{t_n}^s \phi(\hat{X}_n) \hat{X}_rdr + \int_{t_n}^s \sigma K \hat{X}_r d W_r\\
& = & \hat{X}_{n}\exp\{(\phi(\hat{X}_n)-\frac{1}{2}\sigma^2 K^2)(s-t_n) + \sigma K(W_s - W_n)\},
\eeao
\cf \cite[Section 4.4]{kloeden_platen:1995}  which suggests the scheme
\beqq\label{eq:SDscheme}
\hat{X}_{n+1} =  \hat{X}_{n}\exp\{(\phi(\hat{X}_n)-\frac{1}{2}\sigma^2 K^2)\Del + \sigma K\Del W_n\},
\eeqq
where $X_n:= X_{t_n}.$ 

We can show that the numerical scheme also possesses finite moment bounds of any order, that is
\beqq \label{eq:MomentboundSD}
\bfE (\hat{X}_{n})^p\leq \wt{C}_p,
\eeqq
where
$$
\wt{C}_p := \hat{X}_{0}^p \exp\left\{\left(\eta + \frac{\sigma^2K^2}{2}\right)Tp + \frac{p^2\sigma^2K^2}{2}T\right\}
$$
Indeed, by considering $\hat{X}_{0}$ nonnegative then all $\hat{X}_{n}$ are nonnegative and we can write
\beao
\hat{X}_{n+1} &=&  \hat{X}_{n}\exp\left\{\left(\eta  + (\eta - \beta)\hat{X}_n + \frac{\sigma^2 K^2 \hat{X}_n}{1 + \hat{X}_n} -\frac{1}{2}\sigma^2 K^2\right)\Del + \sigma K\Del W_n\right\}\\
&\leq& \hat{X}_{n}\exp\left\{\left(\eta  + \frac{\sigma^2 K^2}{2}\right)\Del\right\}\exp\{\sigma K\Del W_n\}\\
&\leq& \hat{X}_{n-1}\exp\left\{\left(\eta  + \frac{\sigma^2 K^2}{2}\right)2\Del\right\}\exp\{\sigma K(\Del W_n + \Del W_{n-1})\}\\
&\leq& \hat{X}_{0}\exp\left\{\left(\eta  + \frac{\sigma^2 K^2}{2}\right)(n+1)\Del\right\}\exp\{\sigma K \sum_{i=0}^{n+1}\Del W_i\}
\eeao
by using repeatedly (\ref{eq:SDscheme}). Raising to the power of $p$ and taking expectations in the above inequality we reach  (\ref{eq:MomentboundSD}) where we have used the fact the exponential moments of a Wiener process are bounded, that is $\bfE e^{cW_T}\leq e^{c^2T/2},$ for $c\in\bbR.$

Having followed the semi-discrete method we immediately get  a strong convergence result of type 
for the transformed process, that is (see \cite[Theorem 2.1]{halidias_stamatiou:2016}) 
\beqq\label{eq:strongconvSIS}
\bfE \sup_{t\in[0,T]}|\hat{z}_t - \hat{X}_t|^2\rightarrow0 \quad\mbox{for}\quad\Del\downarrow0.
\eeqq

Now we want to reveal the order of convergence in  (\ref{eq:strongconvSIS}), that is find the value of $r$ such that
\beqq\label{eq:strongconvorderSIStrans}
\bfE \sup_{t\in[0,T]}|\hat{z}_t - \hat{X}_t|^2\leq C\Del^{2r}.
\eeqq

Then, by the mean value theorem we have that 
$$
\frac{\hat{z}_t}{1 + \hat{z}_t} - \frac{\hat{X}_t}{1 + \hat{X}_t} = \frac{1}{(1 + v_t)^2}(\hat{z}_t - \hat{X}_t) 
$$
where $v_t$ is between $\hat{z}_t$ and $\hat{X}_t$ so that 
\beam\nonumber  
\bfE \sup_{t\in[0,T]}|I_t - \hat{Y}_t|^2  & = &  \left|\frac{K\hat{z}_t}{1 + \hat{z}_t} - \frac{K\hat{X}_t}{1 + \hat{X}_t}\right|^2\\
\nonumber& \leq & K^2\bfE \sup_{t\in[0,T]}|\hat{z}_t - \hat{X}_t|^2\\
\label{eq:strongconvorderSIS}&\leq& C\Del^{2r}.
\eeam

We show that (\ref{eq:strongconvorderSIS}) holds with $r=1.$ We present our first result.

\bth\label{th1}
For any $q>0,$ the numerical scheme (\ref{eq:SIS_SD}) strongly converges, in $\bbl_q$ sense, to the solution of (\ref{eq3}) with order $1$, that is 
\beqq\label{eq:qconvergence}
\bfE \sup_{n=0,1,\ldots,\ceil{T/\Del}} |\hat{Y}_{t_n} - I_{t_n}|^{q}  \leq C\Del^{q}.
\eeqq
\ethe

Furthermore, we examine the stability behavior of the method. Recall model (\ref{eq3}). We denote the deterministic reproduction number by $\bbr_0^D = \frac{\beta}{K(b+\ga)}.$

\bth[see Theorem 4.1 in \cite{gray_elal:2011}]\label{th:SISlimitStoch}
Let
$$
\bbr_0^S:=\bbr_0^D - \frac{\sigma^2K^2}{2(b+\ga)}<1 \quad \mbox{ and } \quad \sigma^2\leq \frac{\be}{K^2}
$$
Then for any $I_0\in(0,K)$ the solution process of (\ref{eq3}) tends exponentially fast to zero a.s., that is
\beqq\label{eq:exponentialstab}
\limsup_{\tto} \frac{1}{t}\ln(I_t) \leq \be - b - \ga - \frac{1}{2}\sigma^2K^2<0
\eeqq
\ethe
Relation (\ref{eq:exponentialstab}) shows that the disease will die out with probability $1.$ We show in the next result that the proposed numerical scheme inherits perfectly this property, in the following sense. 

\bth\label{th2}
Let the assumptions of Theorem~\ref{th:SISlimitStoch} hold, where also $\sigma^2K^2\leq(b+\ga).$  Then for any $0<\Del<1$ it holds that 
\beqq\label{eq:exponentialstabSD}
\limsup_{n\Del\rightarrow\infty} \frac{1}{n\Del}\ln(\hat{Y}_n) \leq \be - b - \ga - \frac{1}{2}\sigma^2K^2<0
\eeqq
that is the approximation process (\ref{eq:SIS_SD}) tends exponentially fast to zero a.s.
\ethe

In other words the proposed numerical scheme reproduces in a perfect way the stability property of the solution process, when this extra condition on the parameters holds.

\section{Strong Convergence of the method}\label{sec:Con}

Recall (\ref{eq:SIS_SD}). We will compare the proposed scheme with scheme (\ref{eq:SISGrayYang}) which has order $1.$ We are interested in the estimation of 
\beqq\label{eq:closeness of schemes}
|u_n|^q := |\hat{Y}_{n} - Y_n|^q = K^q \left|\frac{\hat{X}_n}{1 + \hat{X}_n} - \frac{e^{X_n}}{1 + e^{X_n}}\right|^q
\eeqq
for any $q>0$ (or at least for $q = 2$), where $\hat{X}_{n}$ is given by (\ref{eq:SDscheme}) and 
$X_{n}$ by (\ref{eq:SISGrayYangTransf}). Actually by \cite[Theorem 2.1]{YANG:2021} 
$$
\bfE \sup_{n=0,1,\ldots,\ceil{T/\Del}} |I_{t_n} - Y_{t_n}|^{q}  \leq C\Del^{q},
$$
therefore we would like to show a result of the following type
\beqq\label{eq:closeness of schemesSUP}
\bfE \sup_{n=0,1,\ldots,\ceil{T/\Del}} |\hat{Y}_{t_n} - Y_{t_n}|^{q}  \leq C\Del^{q},
\eeqq
since then by the triangle inequality we would have
$$
\bfE \sup_{n=0,1,\ldots,\ceil{T/\Del}} |I_{t_n} - \hat{Y}_{t_n}|^{q}  \leq C\Del^{q}.
$$

First note that 
$$\left|\frac{\hat{X}_n}{1 + \hat{X}_n} - \frac{e^{X_n}}{1 + e^{X_n}}\right| \leq 1$$
and if we write $\hat{X}_{n} = e^{\Theta_n},$ with
\beqq\label{eq:log}
\Theta_{n} = \ln(\hat{X}_{n-1}) + (\phi(\hat{X}_{n-1})-\frac{1}{2}\sigma^2 K^2)\Del + \sigma K\Del W_{n-1}
\eeqq
then by application of the mean value theorem for the exponential function we get
$$|u_n|  = K \frac{|e^{\Theta_n} - e^{X_n}|}{(1 + \hat{X}_n)(1 + e^{X_n})}  \leq K |\Theta_n - X_n|.$$
Moreover
\beao
\Theta_n - X_n & = & \ln(\hat{X}_{n-1}) + \left(\eta  + (\eta - \beta)\hat{X}_{n-1} + \frac{\sigma^2 K^2 \hat{X}_{n-1}}{1 + \hat{X}_{n-1}} -\frac{1}{2}\sigma^2 K^2\right)\Del + \sigma K\Del W_{n-1}\\
&& - X_{n-1} - \left(\eta - (b + \gamma)e^{X_{n-1}} + \frac{\sigma^2 K^2}{2} - \frac{\sigma^2 K^2}{1 + e^{X_{n-1}}}\right)\Del - \sigma K\Del W_{n-1}\\
& = & \ln(\hat{X}_{n-1}) - X_{n-1} - (b + \gamma)(\hat{X}_{n-1}  - e^{X_{n-1}})\Del - \sigma^2 K^2\Del\\
&& + \sigma^2 K^2\left( \frac{ \hat{X}_{n-1}}{1 + \hat{X}_{n-1}}  + \frac{1}{1 + e^{X_{n-1}}}\right)\Del\\ 
& = &  \ln(\hat{X}_{n-1}) - X_{n-1} - (b + \gamma)(e^{\Theta_{n-1}}  - e^{X_{n-1}})\Del \\
&&+ \sigma^2 K^2\left( \frac{ \hat{X}_{n-1}}{1 + \hat{X}_{n-1}}  + \frac{1}{1 + e^{X_{n-1}}} - 1\right)\Del\\
& = &  \ln(\hat{X}_{n-1}) - X_{n-1} - (b + \gamma)(e^{\Theta_{n-1}}  - e^{X_{n-1}})\Del + \sigma^2 K^2\left( \frac{ \hat{X}_{n-1}}{1 + \hat{X}_{n-1}}  - \frac{e^{X_{n-1}}}{1 + e^{X_{n-1}}}\right)\Del 
\eeao
where $\eta - \beta = - (b + \gamma)$ or in terms of $v_n := \Theta_n - X_n$
\beqq
v_n  = v_{n-1} - (b + \gamma)(e^{\Theta_{n-1}}  - e^{X_{n-1}})\Del + \sigma^2 K^2u_{n-1}\Del 
\eeqq

Taking the square of the each side of the above equality yields
\beao
(v_n)^2 &=&  (v_{n-1})^2 +  (b + \gamma)^2(e^{\Theta_{n-1}}  - e^{X_{n-1}})^2\Del^2  + \sigma^4 K^4(u_{n-1})^2\Del^2\\
&& + 2\sigma^2 K^2v_{n-1}u_{n-1}\Del - 2(b + \gamma)v_{n-1} (e^{\Theta_{n-1}}  - e^{X_{n-1}})\Del \\
&&- 2(b + \gamma)\sigma^2 K^2(e^{\Theta_{n-1}}  - e^{X_{n-1}})u_{n-1}\Del^2\\
&\leq& (1 + \sigma^4 K^6\Del^2 + 2\sigma^2K^3\Del)(v_{n-1})^2 
+ (b + \gamma)^2e^{2\xi_n}(v_{n-1})^2\Del^2
\\
&\leq& (1 + C\Del)(v_{n-1})^2  + C\Del^2e^{2\xi_n}(v_{n-1})^2,
\eeao
where we have used that $|u_n| \leq K|v_n|$ and the mean value theorem for the exponential function to remove the last two negative terms, $\xi_n$ is between $\Theta_{n-1}$ and $X_{n-1}.$ 
The inequality for $(v_n)^2$ becomes
\beao
(v_n)^2 &\leq& (1 + C\Del)(v_{n-1})^2  + C\Del^2e^{2\xi_n}(v_{n-1})^2\\
&\leq& C\Del^2\sum_{j=0}^{n-1}( 1 + C\Del)^{n-j-1} e^{2\xi_j}(v_{j})^2
\eeao
where  $\xi_j$ is between $\Theta_{j}$ and $X_{j}$. We note that $\sup_{\Del\in(0, 1)}\sup_{n=1,2,\ldots\ceil{T/\Del}}(1 + C\Del)^{n}<\infty.$ Now, we  raise to the power of $q,$ with $q\geq1,$ take the supremum over all $n = 0,1,\ldots,k$ where $k=0,1,\ldots,\ceil{T/\Del},$ and then expectation to the above inequality to find 

\beao
\bfE\sup_{n=0,1,\ldots,k}(v_{n})^{2q} &\leq& C\Del^{2q}\bfE\sup_{n=0,1,\ldots,k}\left| \sum_{j=0}^{n-1}e^{2\xi_j} (v_{j})^2\right|^q\\
&\leq& C\Del^{2q}\Del^{1- q}\bfE\left(\sum_{j=0}^{k-1}e^{2q\xi_j}  (v_{j})^{2q}\right)\\
&\leq& C\Del^{1 + q}\sum_{j=0}^{k-1}\bfE \left(e^{2q\xi_j}  (v_{j})^{2q}\right).
\eeao
We bound the term inside the expectation in the following way
$$
e^{2\xi_j}  (v_{j})^{2} = e^{2\xi_j}(1 + e^{\Theta_j})(1 + e^{X_j})\frac{|u_j|}{K} |v_j| \leq C|v_j| + Ce^{2(\Theta_j\vee X_j)} |v_j|,  
$$
where $a\vee b = \max\{a,b\}.$ Applying  Holder's inequality and Young's inequality we get the following bound for the expectation inside the sum
\beao
\bfE \left(e^{2q\xi_j}  (v_{j})^{2q}\right) &\leq& C\bfE (|v_{j}|^{q}) + C\bfE \left(e^{2q(\Theta_j\vee X_j)}|v_{j})|^{q}\right)\\
&\leq &C\bfE (|v_{j}|^{q}) + C\sqrt{\Del^{q}\bfE \left(e^{4q(\Theta_j\vee X_j)}\right)}\sqrt{\Del^{-q}\bfE (v_{j})^{2q}}\\
&\leq & C\bfE (|v_{j}|^{q}) + C\Del^q + C\Del^{-q} \bfE (v_{j})^{2q}.
\eeao

Collecting all the above estimates we conclude  that
\beao
\bfE\sup_{n=0,1,\ldots,k}(v_{n})^{2q} &\leq& C\Del\sum_{j=0}^{k-1}\bfE(v_{j})^{2q}  +   C\Del^{2q}\\
&\leq& C\Del^{2q}e^{\sum_{j=0}^{k-1}C\Del}\leq C\Del^{2q},
\eeao
where in  the last step we have applied the discrete version of the Gronwall inequality. Thus using once more  $|u_n| \leq K|v_n|$ we find 
$$
\bfE \sup_{n=0,1,\ldots,\ceil{T/\Del}} |u_n|^{2q}  \leq C\Del^{2q},
$$
or in other words (\ref{eq:closeness of schemesSUP}). Therefore (\ref{eq:qconvergence}) is true.

\section{Stability of the method}\label{sec:Stab}

Recall representation (\ref{eq:log}), which we rewrite as
$$
\ln(\hat{X}_{n}) = \ln(\hat{X}_{n-1}) + \left(\eta  - (b+\ga)\hat{X}_{n-1} + \frac{\sigma^2 K^2 \hat{X}_{n-1}}{1 + \hat{X}_{n-1}} -\frac{1}{2}\sigma^2 K^2\right)\Del + \sigma K\Del W_{n-1}.
$$
We will work with $\ln(\hat{X}_{n})$ reaching a result of the type
\beqq\label{eq:lnX}
\limsup_{n\Del\rightarrow\infty} \frac{1}{n\Del}\ln(\hat{X}_n) \leq \be - b - \ga - \frac{1}{2}\sigma^2K^2<0.
\eeqq

After that, the desired inequality is valid since by writing
\beao
\ln(\hat{Y}_n) &=& \ln \left(K\frac{\hat{X}_n}{1 + \hat{X}_n}\right)\\
&\leq& \ln K + \ln \hat{X}_n
\eeao
we immediately get 
\beao
\limsup_{n\Del\rightarrow\infty} \frac{1}{n\Del}\ln(\hat{Y}_n) &\leq& \limsup_{n\Del\rightarrow\infty} \frac{1}{n\Del}\ln K + \limsup_{n\Del\rightarrow\infty} \frac{1}{n\Del}\ln(\hat{X}_n)\\
&\leq& \be - b - \ga - \frac{1}{2}\sigma^2K^2<0.
\eeao

\subsection*{Proof of (\ref{eq:lnX})} We bound $\ln(\hat{X}_{n})$ in the following way

\beao
\ln(\hat{X}_{n}) &\leq & \ln(\hat{X}_{n-1}) + \left(\eta  -\frac{1}{2}\sigma^2 K^2\right)\Del + (\sigma^2K^2-(b+\ga))\hat{X}_{n-1}\Del + \sigma K\Del W_{n-1}\\
&\leq & \ln(\hat{X}_{n-2}) + \ln(\hat{X}_{n-1}) + \left(\eta  -\frac{1}{2}\sigma^2 K^2\right)2\Del + \sigma K\Del W_{n-2}  + \sigma K\Del W_{n-1}\\
&\leq & \ln(\hat{X}_0) + n\Del\left(\eta  -\frac{1}{2}\sigma^2 K^2\right) + \sigma K\sum_{i=0}^{n-1}\Del W_{i},
\eeao
where we used that $\sigma^2K^2\leq b + \ga.$ Therefore
\beqq\label{eqXn-}
\limsup_{n\Del\rightarrow\infty} \frac{1}{n\Del}\ln(\hat{X}_{n}) \leq \left(\eta  -\frac{1}{2}\sigma^2 K^2\right) + \sigma K\limsup_{n\Del\rightarrow\infty} \frac{1}{n\Del}\sum_{i=0}^{n-1}\Del W_{i}
\eeqq

At this point we will use properties of the Wiener process. The sum in the last term of (\ref{eqXn-}) behaves like $W_{n\Del}.$ By the law of the iterated logarithm  \cf \cite[Theorem 1.4.2]{mao:2007}
$$
\limsup_{\tto} \frac{W_t}{\sqrt{2t\ln\ln(t)}} = 1 \quad\mbox{a.s.}
$$
Thus (\ref{eqXn-}) is further bounded as
\beao
\limsup_{n\Del\rightarrow\infty} \frac{1}{n\Del}\ln(\hat{X}_{n}) &\leq& \eta  -\frac{1}{2}\sigma^2 K^2 + \sqrt{2}\sigma K\limsup_{n\Del\rightarrow\infty} \frac{W_{n\Del}}{\sqrt{ 2n\Del\ln\ln(n\Del)}}\limsup_{n\Del\rightarrow\infty} \frac{\sqrt{ \ln\ln(n\Del)}}{\sqrt{n\Del}}\\
&\leq& \eta  -\frac{1}{2}\sigma^2 K^2 + \sqrt{2}\sigma K\limsup_{n\Del\rightarrow\infty}\sqrt{\frac{ \ln\ln(n\Del)}{n\Del}}\\
&\leq& \eta  -\frac{1}{2}\sigma^2 K^2. 
\eeao

\section{Numerical Experiment}\label{sec:Num}

\begin{figure}[ht]
	\centering
	\includegraphics[width=1\textwidth]{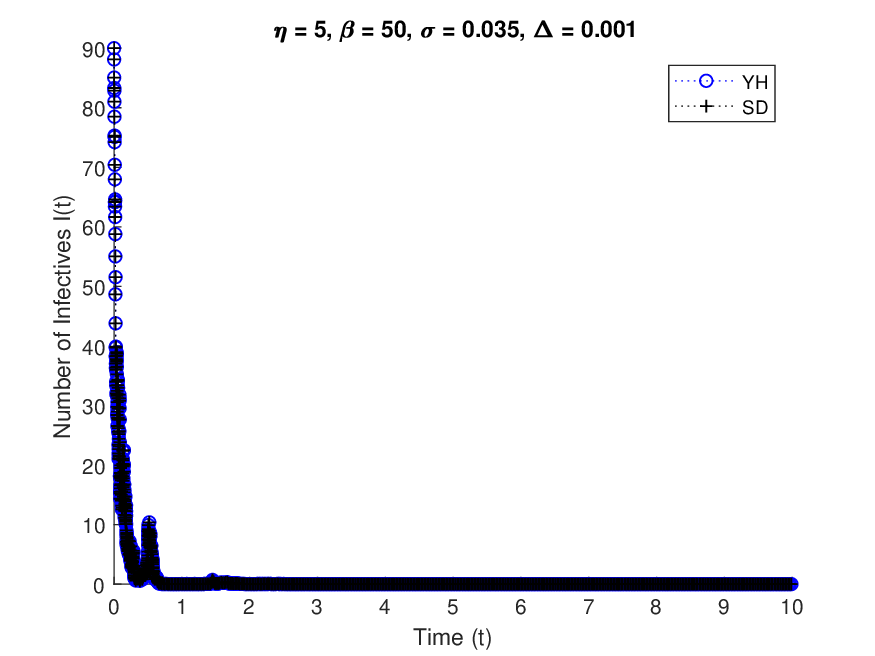}
	\caption{A sample path produced by (\ref{eq:SISGrayYang}) and (\ref{eq:SIS_SD}) using the same Wiener Process for the approximation of model (\ref{eq2}) with parameters as in \cite[Example 3.1]{YANG:2021} and step-size $\Del = 0.001$.}\label{fig:SIS_Yang_SD}
\end{figure}

\begin{figure}[ht]
	\centering
	\includegraphics[width=1\textwidth]{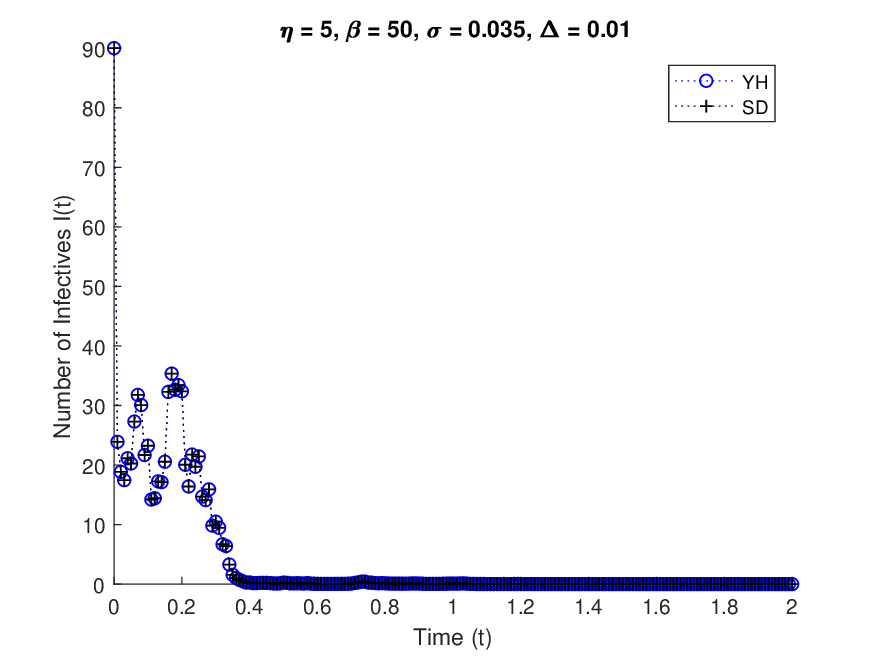}
	\caption{A sample path produced by (\ref{eq:SISGrayYang}) and (\ref{eq:SIS_SD}) using the same Wiener Process for the approximation of model (\ref{eq2}) with parameters as in \cite[Example 3.1]{YANG:2021} and step-size $\Del = 0.01$.}\label{fig:SIS_Yang_SD2}
\end{figure}

We work with an example considered in \cite[Example 3.1]{YANG:2021} and compare the two schemes, see Figures~\ref{fig:SIS_Yang_SD} and ~\ref{fig:SIS_Yang_SD2}.

\begin{figure}[ht]
	\centering	\includegraphics[width=1\textwidth]{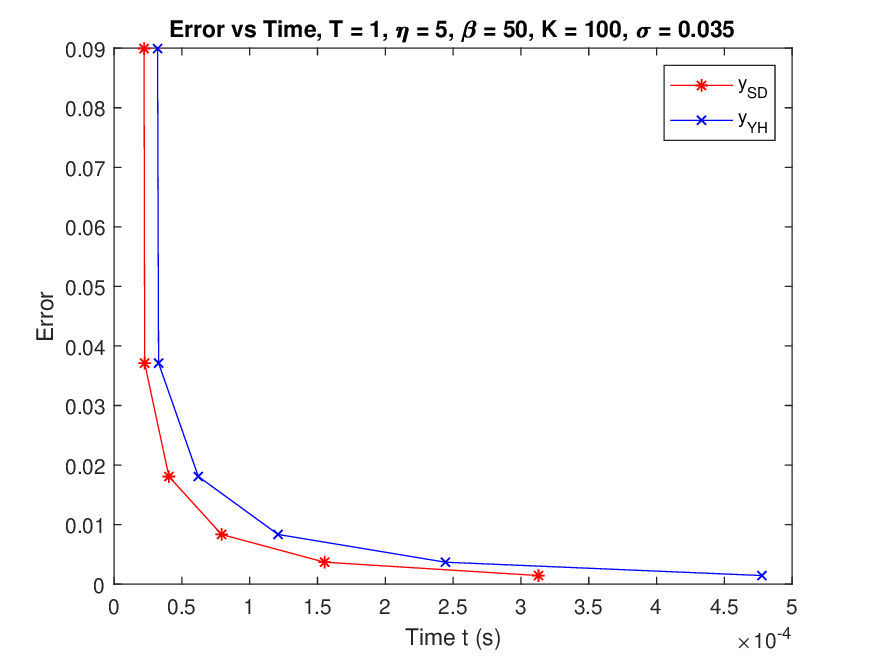}
	\caption{Error as a function of CPU time (in sec) using schemes (\ref{eq:SIS_SD}) and (\ref{eq:SISGrayYang}) for the approximation of (\ref{eq2}).}\label{errortime_SD_YH}
\end{figure}

The proposed domain preserving scheme, seems to produce paths ``close" to the ones found \cite[Example 3.1]{YANG:2021}.  Therefore, apart from an alternative option to qualitative approximation of the solution process, the proposed scheme seems to be superior with respect to computational time, see Figure~\ref{errortime_SD_YH}.
We refer also to the closely related work in \cite{chen_etal:2021} where an explicit truncated method is also used to approximate (\ref{eq:SIStransformed}).

\begin{figure}[ht]
	\centering
	\includegraphics[width=1\textwidth]{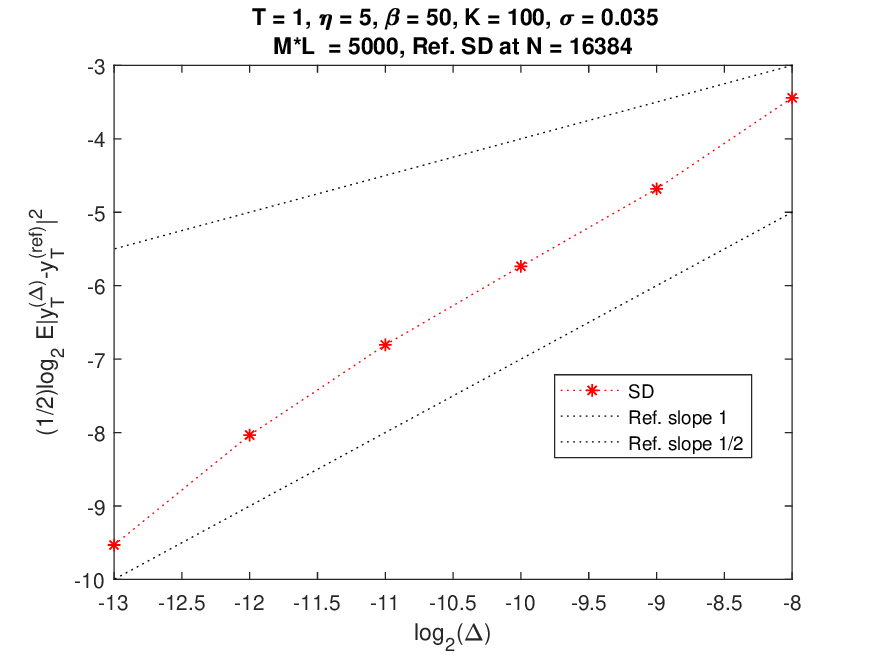}
	\caption{Experimental Order of Convergence for SIS model (\ref{eq2}) for scheme (\ref{eq:SIS_SD}).}\label{fig:SISorderSD}
\end{figure}

\begin{figure}[ht]
	\centering
	\includegraphics[width=1\textwidth]{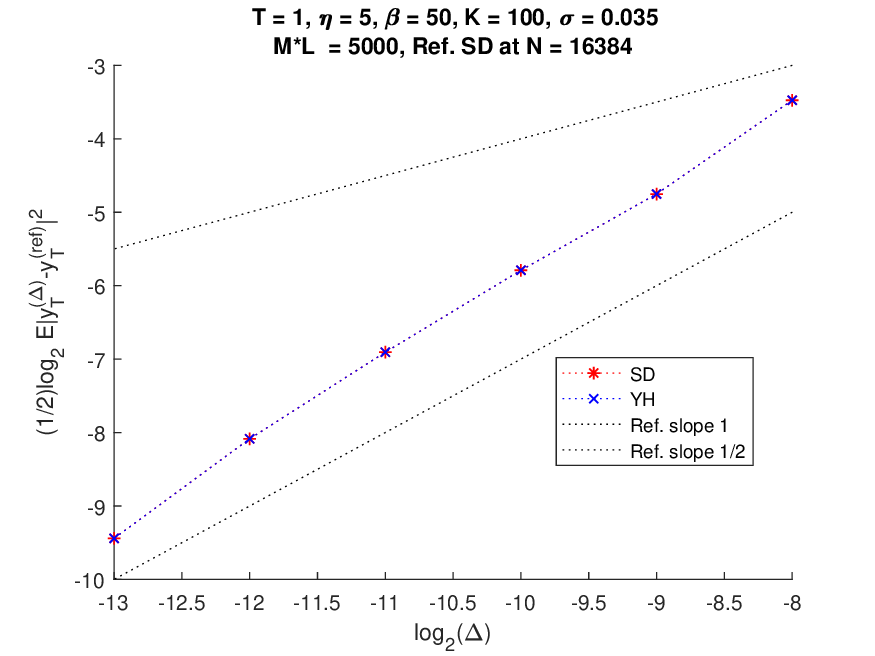}
	\caption{Experimental Order of Convergence for SIS model (\ref{eq2}) for schemes (\ref{eq:SIS_SD}) and (\ref{eq:SISGrayYang}) with YH as reference solution.}\label{fig:SISSD_YH_order}
\end{figure}

\begin{figure}[ht]
	\centering
	\includegraphics[width=1\textwidth]{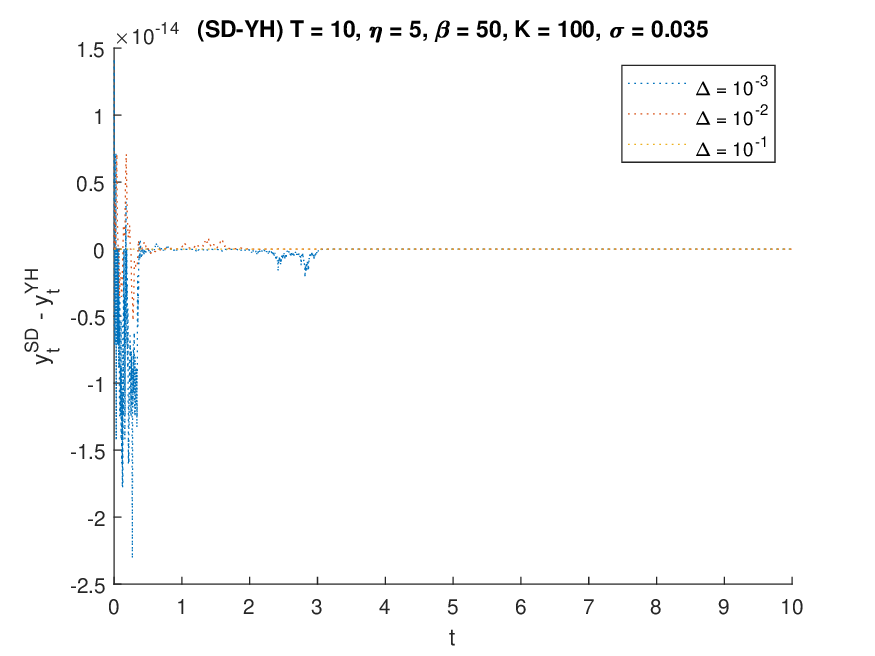}
	\caption{Difference between (\ref{eq:SIS_SD}) and (\ref{eq:SISGrayYang}) for the approximation of (\ref{eq2}).}\label{fig:SISSDminYH}
\end{figure}

\section*{Acknowledgment}
  \noindent The authors wish to acknowledge fruitful discussions with
  A and B.




\bibliographystyle{plain}
\bibliography{sample.bib}


\end{document}